\documentclass[12pt]{amsart}
\usepackage{latexsym,amsmath,amsthm,verbatim,ifthen,amssymb}
\usepackage[latin1]{inputenc}
\usepackage[all]{xy}
\usepackage{graphicx}
\usepackage{epsfig}
\newdir{ >}{{}*!/-7pt/\dir{>}}
\newtheorem{theorem}{Theorem}
\newtheorem{corollary}[theorem]{Corollary}
\newtheorem{lemma}[theorem]{Lemma}
\newtheorem{proposition}[theorem]{Proposition}

\newtheorem{definition}[theorem]{Definition}
\newtheorem{remark}[theorem]{Remark}

\def\secat{{\rm{secat}\hskip1pt}}

\newcommand{\eqr}[1]{\rotatebox{#1}{{\tiny $\!\! \sim \!\!\!$}}}

\def\Q{{\mathbb{Q}}}

\begin{document}

\title{Abstract Sectional Category}
\author[F.J. D\'{\i}az, J.M. Calcines, P.R. D\'{\i}az, A. Murillo Mas and J. Remedios]
{F. D\'{\i}az D\'{\i}az, J.M. Garc\'{\i}a Calcines, \\ P.R.
Garc\'{\i}a D\'{\i}az, A. Murillo Mas and J. Remedios G\'{o}mez}



\address{F.J. D\'{\i}az; J.M. Garc\'{\i}a Calcines; \newline \indent
P.R. Garc\'{\i}a D\'{\i}az; J. Remedios G\'omez \newline \indent
Departamento de Matem\'{a}tica Fundamental \newline \indent
Universidad de La Laguna \newline \indent 38271 La Laguna, Spain.}
\email{fradiaz@ull.es} \email{jmgarcal@ull.es}
\email{prgdiaz@ull.es} \email{jremed@ull.es}

\address{A. Murillo Mas \newline
\indent Departamento de Algebra, Geometr\'{\i}a  y Topolog\'{\i}a
\newline \indent Universidad de M\'alaga \newline \indent Ap.59,
29080 M\'alaga, Spain.} \email{aniceto@agt.cie.uma.es}

\thanks{Partially supported by FEDER, the Ministerio de Educaci\'on y
Ciencia grants MTM2009-12081, MTM2010-18089 and by the Junta de
Andaluc\'\i a grant FQM-213}

\begin{abstract}
We study, in an abstract axiomatic setting, the notion of
sectional category of a morphism. From this, we unify and
generalize known results about this invariant in different
settings as well as we deduce new applications.
\end{abstract}

\maketitle

\section*{Introduction}
The sectional category $\mbox{secat}\hspace{2pt}(p)$ of a
fibration $p:E\twoheadrightarrow B$, originally introduced by A.
Schwarz \cite{Sch}, is defined as the least integer $n$ such that
$B$ admits a cover constituted by $n+1$ open subsets, on each of
which $p$ has a local section. It is a lower bound of the
Lusternik-Schnirelmann category of the base space and it is also a
generalization of this invariant since
$\mbox{secat}\hspace{2pt}(p)=\mbox{cat}\hspace{2pt}(B)$ when $E$
is contractible.  Apart from the original applications of the
sectional category in the  classification of bundles or the
embedding problem \cite{Sch}, this numerical invariant has proved
to be useful in different settings. For instance, Smale \cite{Sm}
showed that the sectional category of a certain fibration provides
a lower bound for the complexity of algorithms computing the roots
of a complex polynomial. We can also mention the work of M. Farber
\cite{Far,Far2} who introduced the topological complexity of a
given space $X$ as the sectional category of the path fibration
$X^I\rightarrow X\times X,$ $\alpha \mapsto (\alpha (0),\alpha
(1)).$ In robotics, when $X$ is thought to be the configuration
space associated to the motion of a given mechanical system, this
invariant measures, roughly speaking, the minimum amount of
instructions of any algorithm controlling the given system.

In general, the sectional category of a fibration is hard to
compute. The notion of Lusternik-Schnirelmann category (L.-S.
category, for short) has the same disadvantage. In order to face
this problem for L.-S. category there have been several attempts
to describe it in a more functorial and therefore manageable form;
among the most successful ones we can mention the Whitehead and
Ganea characterizations. Many other approximations of L.-S.
category have been introduced. One of them relies in an important
algebraic technique for obtaining lower bounds. It consists of
taking models of spaces in an algebraic category where a notion of
L.-S.-category type invariant is given. Such algebraic category
must posses an abstract notion of homotopy, usually established in
an axiomatic homotopy setting, such as a Quillen model category.
Then the algebraic L.-S. category of the model of X is a lower
bound of the original L.-S. category of $X.$ During the progress
of this technique, several algebraic notions of L.-S. category
have been appearing. In 1993, in order to give a common point for
all of them, Doeraene \cite{D} introduced the notion of L.-S.
category in a Quillen model category. Actually, in his work
Doeraene develops two different notions of L.-S. category, which
are the analogous to the Ganea and Whitehead characterizations in
the topological case and proves that, under the crucial cube
axiom, these notions agree, as expected. As far as the sectional
category is concerned, not much has been done in this direction.
In the work of A. Schwarz \cite{Sch} it was established a
Ganea-type characterization of sectional category. Namely, if
$p:E\twoheadrightarrow B$ is a fibration we can consider
$j_n:*^n_BE\rightarrow B,$ which is the $n$-th fold join of $p.$
If the base space $B$ is paracompact, then A. Schwarz proved that
$\mbox{secat}\hspace{2pt}(p)\leq n$ if and only if $j_n$ admits a
(homotopy) section. Clapp and Puppe \cite[Cor. 4.9]{clpp} also
obtained a Whitehead-type characterization of sectional category;
more precisely, for a given map $p:E\rightarrow B$ with associated
cofibration $\hat{p}:E\rightarrow \hat{B},$
$\mbox{secat}\hspace{2pt}(p)\leq n$ if and only if the diagonal
map $\Delta _{n+1}:\hat{B}\rightarrow \hat{B}^{n+1}$ factors, up
to homotopy, through the $n$-th fat wedge
$T^n(\hat{p})=\{(b_0,b_1,...,b_n)\in \hat{B}^{n+1}:x_i\in
E,\hspace{3pt}\mbox{for some}\hspace{3pt}i \}.$ With this
characterization Fass\`o \cite{fa} studied the sectional category
of the corresponding algebraic model of $p$ in rational homotopy.
These functorial characterizations in the topological case open a
door through an axiomatization of sectional category. In this
direction an initial advance has been made by T. Kahl in \cite{K}.
In his work he gives the notion of abstract sectional category
through a certain variation of inductive L.-S. category in the
sense of Hess-Lemaire \cite{He-L}.

Our aim in this paper  is to develop, in the same spirit as
Doeraene did in \cite{D} with the L.-S. category, the notion of
sectional category in an abstract homotopy setting and to deduce
some applications. In the first section we recall some background
to set the axiomatic framework in which we shall work as well as
the main tools that will be used. In \S2 we introduce, under two
different approaches, the concept of sectional category of a given
morphism. Then, in \S3 we present the main properties of this
invariant and finally, in the fourth section, we give some
applications.

\section{Preliminaries: $J$-category and main notions.}

In this paper we shall work in a $J$-category \cite{D}, which
includes the cases of a pointed cofibration and fibration category
in the sense of Baues \cite{B} or a pointed proper model category
\cite{Q,Q2} as long as they satisfy the ``cube lemma". The aim of
this section is to provide some of the most important notions and
properties given in such a homotopy setting. For proofs and more
details the reader is referred to Doeraene's paper \cite{D} or his
thesis \cite{D2}.

Explicitly, a $J$-category $\mathcal{C}$ is a category  with a
zero object $0$ and endowed with three classes of morphisms called
fibrations ($\twoheadrightarrow $), cofibrations ($\rightarrowtail
$) and weak e\-qui\-va\-lences
($\stackrel{\eqr{0}\:\,}{\rightarrow }$), satisfying the following
set of axioms (J1)-(J5) below. Recall that a morphism which is
both a fibration (resp. cofibration) and a weak equivalence is
called {\it trivial fibration} (resp. {\it trivial cofibration}).
An object $B$ is called {\it cofibrant model} if every trivial
fibration $p:E\stackrel{\eqr{0}\,}{\twoheadrightarrow }B$ admits a
section.

\begin{enumerate}
\item[\textbf{(J1)}] Isomorphisms are trivial
cofibrations and also trivial fibrations. Fibrations and
cofibrations are closed by composition. If any two of $f,$ $g,$
$gf$ are weak equivalences, then so is the third.

\item[\textbf{(J2)}] The pullback of a fibration $p:E\twoheadrightarrow B$
and any morphism $f:B'\rightarrow B$
$$\xymatrix{
  {E'} \ar@{>>}[d]^{\overline{p}} \ar[r]^{\overline{f}}
                      & E \ar@{>>}[d]^{p}    \\
  {B'} \ar[r]_{f}     & B               }
$$
\noindent always exists  and  $\overline{p}$ is a fibration.
Moreover, if $f$ (respec.$p$) is a weak e\-qui\-va\-lence, then so
is $\overline{f}$ (respec. $\overline{p}$). The dual assertion is
also required.

\item[\textbf{(J3)}] For any map $f:X\rightarrow Y$ there exist
an $F$-factorization (i.e., $f=p\tau $ where $\tau $ is a weak
equivalence and $p$ is a fibration) and a $C$-factorization (i.e.,
$f=\sigma i,$ where $i$ is a cofibration and $\sigma $ is a weak
equivalence).

\item[\textbf{(J4)}] For any object $X$ in $\mathcal{C},$ there exists a
trivial fibration
$p_X:\overline{X}\stackrel{\eqr{0}\,}{\twoheadrightarrow }X,$ in
which $\overline{X}$ is a cofibrant model. The morphism
$p_X:\overline{X}\stackrel{\eqr{0}\,}{\twoheadrightarrow }X$ is
called \textit{cofibrant replacement} for $X.$

\end{enumerate}

A commutative square
$$\xymatrix{
  {D} \ar[d]^{g'} \ar[r]^{f'}
                      & C \ar[d]^{g}    \\
  {A} \ar[r]_{f}     & B               }
$$ \noindent is said to be a \textit{homotopy pullback} if for some
(equivalently \textit{any})
$F$-factorization of $g$ (equivalently $f$ or both),  the induced map 
from $D$ to the pullback $E'=A\times _B E$ is a weak equivalence
$$\xymatrix@C=.7cm@R=.5cm{
{D} \ar[dd]_{g'} \ar[rrr]^{f'} \ar@{.>}[dr] & & & {C} \ar[dd]^g
\ar[dl]^{\tau }_{\eqr{45}} \\
 & {E'} \ar@{>>}[dl]^{\overline{p}} \ar[r]_{\overline{f}} & {E} \ar@{>>}[dr]^p & \\
 {A} \ar[rrr]_f & & & B  }$$
The notion of \textit{homotopy pushout} is dually defined.

 \begin{enumerate}

 \item[\textbf{(J5)}] {\it The cube axiom}. Given any commutative cube where the
bottom face is a homotopy pushout and the vertical faces are
homotopy pullbacks, then the top face is a homotopy pushout.

\end{enumerate}

\begin{remark}
As pointed out by Doeraene, (J1)-(J4) axioms allow us to replace
'some' by 'any' in the definition of homotopy pullback, or to use
an $F$-factorization of $f$ instead of $g.$
\end{remark}

We are particularly interested in knowledge of objects and
morphisms \textit{up to weak equivalence}. Two objects $A$ and
$A'$ in $\mathcal{C}$ are said to be {\it weakly equivalent} if
there exists a finite chain of weak equivalences joining $A$ and
$A'$
$$\xymatrix{
A \ar@{-}[r]^{\eqr{0}} & \bullet \ar@{-}[r]^(.33){\eqr{0}} &
\bullet \: \cdots \cdots \: {\bullet }\ar@{-}[r]^(.67){\eqr{0}} &
{A'} }$$ where the symbol $\xymatrix{{\bullet } \ar@{-}[r] &
{\bullet}}$ means an arrow with either left or right orientation.
One can analogously define the notion of \textit{weakly equivalent
morphisms} by considering a finite chain of weak equivalences in
the category $\textrm{Pair}(\mathcal{C})$ of morphisms in
$\mathcal{C}$ (\cite{B} Def. II.1.3)
$$\xymatrix{
A \ar[d]_f \ar@{-}[r]^{\sim } & \bullet \ar[d] \ar@{-}[r]^{\sim} &
\bullet  \ar[d] \ar@{-}[r]^{\sim } &
{\bullet}  \ar[d] \ar@{-}[r]^{\sim } & {\bullet}  \ar[d] \ar@{-}[r]^{\sim } & {A'} \ar[d]^{f'}\\
{B} \ar@{-}[r]_{\sim } &  {\bullet } \ar@{-}[r]_{\sim }  &
{\bullet} \ar@{-}[r]_{\sim } & {\bullet} \ar@{-}[r]_{\sim } &
{\bullet} \ar@{-}[r]_{\sim } & {B'}}$$

\begin{definition} {\rm
Given two morphisms $f:A\rightarrow B$ and $g:C\rightarrow B$,
consider any $F$-factorization of $g=p\tau $ and the pullback of
$f$ and $p.$ Let $\overline{f}$ and $\overline{p}$ the base
extensions of $f$ and $p$ respectively. Then, take any
$C$-factorization of $\overline{f}=\sigma i$ and the pushout of
$\overline{p}$ and $i.$ This pushout object is denoted by $A*_B C$
and is called the {\em join} of $A$ and $C$ over $B.$ The dotted
induced map from $A*_B C$ to $B$ is called the \emph{join morphism
}of $f$ and $g.$
$$\xymatrix@R=.5cm@C=.5cm{
{E'} \ar[rr]^{\overline{f}} \ar@{ >->}[dr]_i
\ar@{>>}[ddd]_{\overline{p}} & & {E} \ar@{>>}[ddd]^p & {C}
\ar[l]_{\tau }^{\eqr{0}} \ar[dddl]^g \\
& {Z} \ar[ur]_{\sigma }^{\eqr{40} } \ar[d] & & \\
& {A*_B C} \ar@{.>}[dr] & & \\
{A} \ar@{ >->}[ur] \ar[rr]_f & & B & }$$ }\end{definition} The
object $A*_B C$ and the join map are well defined and they are
symmetrical up to weak equivalence \cite{D,D2}.

An important result that allows us to see that if a property holds
for some $F$-factorization, then it also holds for any
$F$-factorization is the following lemma. Recall from \cite{B}
that in a fibration category a \emph{relative cocylinder} of a
fibration $p:E\twoheadrightarrow B$ is just an $F$-factorization
of the morphism $(id_E,id_E):E\rightarrow E\times _BE,$ where
$E\times _BE$ denotes the pullback of $p$ with itself
$$\xymatrix{
{E} \ar[rr]^{(id_E,id_E)} \ar[dr]_{\sim } & & {E\times _B E} \\
& {Z_p} \ar@{>>}[ur]_{(d_0,d_1)} & }$$ Then, given
$f,g:X\rightarrow E$ such that $pf=pg,$ it is said that \emph{$f$
is homotopic to $g$ relative to $p$} ($f\simeq g$ rel. $p$) if
there exists a morphism $F:X\rightarrow Z_p$ such that $d_0F=f$
and $d_1F=g.$ When $p=0:E\twoheadrightarrow 0$ is the zero
morphism we obtain the notion of non relative homotopy (and write
$f\simeq g$). In this case, the cocylinder of
$0:E\twoheadrightarrow 0$ will be denoted by $E^I.$

\begin{lemma}\label{lifting-lemma}{\rm \cite[II.1.11]{B}
 Consider a commutative diagram of unbroken arrows:
$$\xymatrix{
{D} \ar[d]^{\sim }_{\tau } \ar[r]^g & {E} \ar@{>>}[d]^p \\
{A} \ar[r]_f \ar@{.>}[ur]_l & {B} }$$
\begin{enumerate}
\item[(a)] If $A$ is a cofibrant model, then there is a morphism
$l:A\rightarrow E$ such that $pl=f.$
\item[(b)] If $A$ and $D$ are cofibrant models, then there is a morphism
$l:A\rightarrow E$ for which $pl=f$ and $l\tau \simeq g$ rel. $p.$
Moreover, if $g$ is a weak equivalence, then so is $l.$
\end{enumerate}}
\end{lemma}

We also recall the notion of weak lifting.
\begin{definition}\label{good-defi-weaklift}{\rm
Let $f:A\rightarrow B$ and $g:C\rightarrow B$ be morphisms in
$\mathcal{C}.$  We say that $f$ admits a \textit{weak lifting
along} $g$ if for some $F$-factorization $g=p\tau $ of $g$ and for
some cofibrant replacement $p_A:\overline{A}\stackrel{\sim
}{\rightarrow }A$ of $A$ there exists a commutative diagram
$$\xymatrix@R=.4cm@C=.4cm{
 & & C \ar[dd]^g \ar[dl]^{\tau }_{\eqr{45}} \\
 & E \ar@{>>}[dr]^p & \\
 {\overline{A}} \ar@{.>}[ur]^s \ar[rr]_{fp_A} &  & B } $$

In the particular case $f=id_B$ we say that $g:C\rightarrow B$
admits a \textit{weak section}.}
\end{definition}

This notion does not depend on the choice of the $F$-factorization
nor on the cofibrant replacement. In order to check this fact one
has to use Lemma \ref{lifting-lemma} above and the following
result. The details are left to the reader.

\begin{lemma}\label{lifting-triv-fib}{\rm \cite[II.1.6]{B}
 Let $p:X\stackrel{\sim }{\twoheadrightarrow }Y$ be a trivial
fibration and $f:A\rightarrow Y$ any morphism, with $A$ a
cofibrant model. Then there exists a lift of $f$ with respect to
$p,$ i.e. a morphism $\tilde{f}:A\rightarrow X$ such that
$p\tilde{f}=f$
$$\xymatrix{
 &  {X} \ar@{>>}[d]^{p}_{\sim } \\
 {A} \ar[r]_{f} \ar@{.>}[ur]^{\tilde{f}} & {Y} }$$}
\end{lemma}

Another important notion that will be used in this paper is the
one of weak pullback.

\begin{definition}\label{weak-pullback}
{\rm Let $f:A\rightarrow B$, $f':A'\rightarrow B'$ and
$b:B\rightarrow B'$ be morphisms in $\mathcal{C}.$ It is said that
$A\mbox{-}A'\mbox{-}B'\mbox{-}B$ is a \textit{weak pullback} if
for some $F$-factorization $f'=p\tau $ and some cofibrant
replacement $p_A:\overline{A}\stackrel{\sim }{\twoheadrightarrow
}A$ of $A$ there exists a homotopy pullback
$$\xymatrix{
{\overline{A}} \ar[d]_{fp_A}^{\hspace{15pt}\mbox{h.p.b.}}
\ar@{.>}[r]^x & X \ar@{>>}[dr]_p & {A'} \ar[d]^{f'} \ar[l]_{\tau
}^{\sim } \\ B \ar[rr]_b & & {B'} }$$ }
\end{definition}

\begin{remark}
Any homotopy pullback is a weak pullback. Again, Lemma
\ref{lifting-lemma}, axiom (J4) and Lemma \ref{lifting-triv-fib}
allow us to replace the word 'some' by 'any' in the above
definition. We also have to take into account that the composition
of homotopy pullbacks is a homotopy pullback (in fact there is a
Prism Lemma for homotopy pullbacks \cite[Prop. 1.1]{D}) and that
the weak equivalences in the category $\textrm{Pair}(\mathcal{C})$
of morphisms in $\mathcal{C}$ are homotopy pullbacks.
\end{remark}

\section{Sectional category. Ganea and Whitehead approaches.}
As in Doeraene's work, from now on we will assume that
$\mathcal{C}$ \emph{is a $J$-category in which all objects are
cofibrant models}. Therefore we will take as cofibrant
replacements the corresponding identities. It is important to
remark that in a general $J$-category we will also obtain the same
results. However, the exposition and/or the arguments in this
general case would be affected by unessential technical
complications. So just for the sake of simplicity and comfort we
admit this assumption without lost of generality. Essentially, the
key point for the pass from our assumption to the general case is
established by considering cofibrant replacements:

\begin{itemize}
\item
\emph{ Any object $X$ in $\mathcal{C}$ has a cofibrant
replacement, that is, a trivial fibration
$p_X:\overline{X}\stackrel{\eqr{0}\,}{\twoheadrightarrow }X,$ in
which $\overline{X}$ is a cofibrant model.} ((J4) axiom )

\item \emph{Any morphism $f:X\rightarrow Y$ in $\mathcal{C}$
has a cofibrant replacement, that is, given cofibrant replacements
$p_X,$ $p_Y$ of $X$ and $Y,$ there exists an induced morphism
$\overline{f}:\overline{X}\rightarrow \overline{Y}$ making
commutative the following square
$$\xymatrix{
{\overline{X}} \ar@{>>}[d]^{\sim }_{p_X} \ar[r]^{\overline{f}} &
{\overline{Y}} \ar@{>>}[d]_{\sim }^{p_Y} \\
{X} \ar[r]_f & {Y} }$$} \end{itemize} Observe that the second item
holds thanks to Lemma \ref{lifting-triv-fib}. Using these simple
facts when necessary and working a little bit harder the reader
should be able to prove our results when not all objects are
cofibrant models.

We are now prepared for the definition of sectional category of a
morphism in $\mathcal{C}$ under two different approaches. In the
following definition, only axioms (J1)-(J4) are needed.

\begin{definition}\label{$G$-sectional category}{\rm
Let $p:E\rightarrow B$ be any morphism in $\mathcal{C}$ (not
necessarily a fibration). We consider for each $n$ a morphism
$h_n:*^n_B E\rightarrow B$ inductively as follows:
\begin{enumerate}
\item $h_0=p:E\rightarrow B$ (so $*^0_B E=E$)

\item Assume that $h_{n-1}:*^{n-1}_B E\rightarrow B$ is already
constructed. Then $h_n$ is the join morphism of $p$ and $h_{n-1}:$
$$\xymatrix@R=.5cm@C=.5cm{
{\bullet } \ar[rr] \ar@{ >->}[dr] \ar@{>>}[ddd] & & {E'}
\ar@{>>}[ddd] & {E}
\ar[l]^{\eqr{0}} \ar[dddl]^p \\
& {\bullet } \ar[ur]^{\eqr{40} } \ar[d] & & \\
& {*^n _B E} \ar@{.>}[dr]^{h_n} & & \\
{*^{n-1}_B E} \ar@{ >->}[ur] \ar[rr]_{h_{n-1}} & & B & }$$
\end{enumerate}

Then, the {\em Ganea sectional category of $p$},
$\mbox{Gsecat}(p)$, is the least integer $n\le \infty$ such that
$h_n$ admits a weak section
$$\xymatrix@R=.4cm@C=.4cm{
 & & {*^n_B E} \ar[dd]^{h_n} \ar[dl]_{\eqr{45}} \\
 & {\bullet } \ar@{>>}[dr] & \\
 {B} \ar@{.>}[ur] \ar@{>>}[rr]_{id_B}^{\sim } &  & B } $$
}
\end{definition}

\begin{remark}
Observe that  $\mbox{Gsecat}(p)=0$ if and only if $p$ has a weak
section. Moreover, in the  topological setting this invariant
coincides with $\mbox{secat}(p),$ the classical sectional category
of a given fibration $p:E\twoheadrightarrow B$, with $B$
paracompact. In fact, the $n$-th iterated join of $p$ over $B$,
$h_n:*^n_B E\rightarrow B$ has a homotopy section if and only if
$B$ can be covered by $n+1$ open subsets, each of them having a
local homotopy section \cite{Ja,Sch}.
\end{remark}

Now we show that this is an invariant up to weak equivalence.

 \begin{proposition}\label{salvoequivdeb}{\rm
If $p:E\rightarrow B$ and $p':E'\rightarrow B'$ are weakly
equivalent morphisms, then $\mbox{Gsecat}(p)=\mbox{Gsecat}(p').$}
\end{proposition}

For the proof we shall use the following result.

\begin{lemma}\label{buenadefi}{\rm \cite[Lemma
3.5]{D} Consider the following commutative diagram in
$\mathcal{C}$
$$\xymatrix{
{A} \ar[d]_x \ar[rr]^f & & {B} \ar[dd]^b & & {C} \ar[ll]_g
\ar[d]^y \\
{\bullet } \ar@{>>}[drr] & & & & {\bullet } \ar@{>>}[dll] \\
{A'} \ar[u]^{\sim } \ar[rr]_{f'} & & {B'} & & {C'} \ar[ll]^{g'}
\ar[u]_{\sim } }$$ That is, $bf$ admits a weak lifting along $f'$
and $bg$ admits a weak lifting along $g'.$ Let $j:A*_BC\rightarrow
B$ and $j':A'*_{B'}C'\rightarrow B'$ denote the corresponding join
maps. Then $bj$ admits a weak lifting along $j'$
$$\xymatrix{
{A*_B C} \ar[d]_{j} \ar[r] & {\bullet }
\ar@{>>}[dr] & {A'*_{B'}C'} \ar[d]^{j'} \ar[l]^(.6){\sim } \\
B \ar[rr]_b & & {B'} }$$ Furthermore, if $b,x$ and $y$ are weak
equivalences, then $A*_B C$ is weakly equivalent to $A'*_{B'}C'$
via the above diagram.}
\end{lemma}


\begin{proof}[Proof of Proposition \ref{salvoequivdeb}]
We can suppose without losing generality that there is a
commutative diagram of the following form
$$\xymatrix{
{E} \ar[d]_p \ar[r]^u_{\sim } & {E'} \ar[d]^{p'} \\
{B} \ar[r]^{\sim }_{v} & {B'} }$$
 Let us see by induction on $n$
that $h_n:*^n_BE\rightarrow B$ and $h'_n:*^n_{B'}E'\rightarrow B'$
are weakly equivalent morphisms. Indeed, for $n=0$ it is certainly
true. Now suppose that $h_{n-1}$ and $h'_{n-1}$ are weakly
equivalent. Again we can assume, without losing generality, that
there is a commutative square
$$\xymatrix{
{*^{n-1}_BE} \ar[d]_{h_{n-1}} \ar[r]^{w}_{\sim } & {*^{n-1}_{B'}E'} \ar[d]^{h'_{n-1}} \\
{B} \ar[r]^{\sim }_{v} & {B'} }$$ Now take
 $h'_{n-1}=q\lambda $ and $p'=r\mu $ $F$-factorizations. Then we
 have a commutative diagram
$$\xymatrix{
{*^{n-1}_BE} \ar[d]_{\lambda w}^{\sim } \ar[rr]^{h_{n-1}} & & {B}
\ar[dd]^v_{\sim } & & {E} \ar[ll]_p
\ar[d]^{\mu u}_{\sim } \\
{\bullet } \ar@{>>}[drr]^q & & & & {\bullet } \ar@{>>}[dll]_r \\
{*^{n-1}_{B'}E'} \ar[u]_{\sim }^{\lambda} \ar[rr]_{h'_{n-1}} & &
{B'} & & {E'} \ar[ll]^{p'} \ar[u]^{\sim }_{\mu } }$$ which,
applying Lemma \ref{buenadefi}, gives rise to this one
$$\xymatrix{
{*^n_BE} \ar[d]_{h_n} \ar[r]_{\sim } & {\bullet }
\ar@{>>}[dr] & {*^n_{B'}E'} \ar[d]^{h'_n} \ar[l]^(.6){\sim } \\
B \ar[rr]_v^{\sim} & & {B'} }$$ \noindent

Now we have that $h_n$ admits a weak section if and only if $h'_n$
admits a weak section. In order to check this assertion, one has
just to take into account Lemma \ref{lifting-lemma} and the fact
that the pullback of $\bullet \twoheadrightarrow B'$ and
$v:B\stackrel{\sim }{\rightarrow }B'$ gives rise to an
$F$-factorization of $h_n$ in a natural way.
\end{proof}

Now we give a Whitehead-type definition of sectional category.

\begin{definition}{\rm
\label{$W$-sectional category}
 Let $p:E\rightarrow B$ be
any morphism in $\mathcal{C}$ where $B$ is $e$-fibrant, that is,
the zero morphism $B\rightarrow 0$ is a fibration. We define
$j_n:T^n(p)\rightarrow B^{n+1}$ inductively as follows:
\begin{enumerate}
\item $j_0=p:E\rightarrow B$ (so $T^0(p)=E$)

\item If $j_{n-1}:T^{n-1}(p)\rightarrow B^n$ is constructed, then
$j_n$ is the following join construction:
$$\xymatrix@R=.5cm@C=.5cm{
{\bullet } \ar[rr] \ar@{ >->}[dr] \ar@{>>}[ddd] & & {\bullet }
\ar@{>>}[ddd] & {B^n\times E}
\ar[l]^(.6){\eqr{0}} \ar[dddl]^{id_{B^n}\times p} \\
& {\bullet } \ar[ur]^{\eqr{40} } \ar[d] & & \\
& {T^n(p)} \ar@{.>}[dr]^{j_n} & & \\
{T^{n-1}(p)\times B} \ar@{ >->}[ur] \ar[rr]_{j_{n-1}\times id_B} &
& B^{n+1} & }$$
\end{enumerate}

Then the {\em Whitehead sectional category of $p$},
$\mbox{Wsecat}(p)$, is the least integer $n\le \infty$ such that
the diagonal morphism $\Delta _{n+1}:B\rightarrow B^{n+1}$ admits
a weak section along $j_n:T^n(p)\rightarrow B^{n+1}:$
$$\xymatrix@R=.4cm@C=.4cm{
 & & {T^n(p)} \ar[dd]^{j_n} \ar[dl]_{\eqr{45}} \\
 & {\bullet } \ar@{>>}[dr] & \\
 B \ar@{.>}[ur] \ar[rr]_{\Delta _{n+1}} &  & {B^{n+1}} } $$ }
\end{definition}

Observe that, in order to define $\mbox{Wsecat}(p),$ we have to
consider $B$ an e-fibrant object to ensure that all products
$B^n,$ $T^n(p)\times B$ and $B^n\times E$ exist ($n\geq 0$). Now
we extend $\mbox{Wsecat}(p)$ to the general case, in which $B$
need not be e-fibrant. For it consider an $F$-factorization
$\xymatrix{B \ar[r]_{\sim }^{\tau } & {F} \ar@{>>}[r] & {0}}$ of
the zero morphism. Then we define
$$\mbox{Wsecat}(p):=\mbox{Wsecat}(\tau p)$$

\begin{lemma}\label{good-defi-whit}
{\rm If $p:E\rightarrow B$ is any morphism, then
$\mbox{Wsecat}(p)$ does not depend on the choice of the
$F$-factorization for $B\rightarrow 0.$}
\end{lemma}

\begin{proof}
Consider $\xymatrix{B \ar[r]_{\sim }^{\tau } & {F} \ar@{>>}[r] &
{0}}$ and $\xymatrix{B \ar[r]_{\sim }^{\tau '} & {F'} \ar@{>>}[r]
& {0}}$ two such $F$-fac\-tor\-izations. Then, by Lemma
\ref{lifting-lemma}(b) applied to the following commutative
diagram
$$\xymatrix{
{B} \ar[r]^{\tau }_{\sim } \ar[d]_{\tau '}^{\sim } & {F}
\ar@{>>}[d] \\ {F'} \ar@{.>}[ur]^{\sim }_h \ar@{>>}[r] & {0} }$$
\noindent there exists a weak equivalence $h:F'\stackrel{\sim
}{\rightarrow }F$ such that $h\tau '\simeq \tau .$ Take a homotopy
$H:B\rightarrow F^I$ verifying that $d_0H=h\tau '$ and $d_1H=\tau
$ and consider the commutative diagram, where the codomain of each
vertical arrow is an e-fibrant object
$$\xymatrix{
 {E} \ar[d]_{\tau ' p} \ar@{=}[r] &
{E} \ar[d]_{h\tau ' p} \ar@{=}[r] & {E} \ar[d]_{Hp} \ar@{=}[r] &
{E} \ar[d]_{\tau p}   \\
{F'}  \ar[r]^{\sim }_h & {F} & {F^I} \ar[l]_{\sim }^{d_0}
\ar[r]^{\sim }_{d_1} & {F} }$$ \noindent This diagram shows that
$\tau p$ and $\tau 'p$ are weakly equivalent morphisms. Observe
that, since $F\times F$ is e-fibrant and by definition there is a
fibration $(d_0,d_1):F^I\twoheadrightarrow F\times F,$ we have
that the cocylinder object $F^I$ is also e-fibrant. Finally,
considering a similar argument to that given in the proof of
Proposition \ref{salvoequivdeb} we obtain the identity
$\mbox{Wsecat}(\tau p)=\mbox{Wsecat}(\tau 'p).$
\end{proof}

\begin{proposition}\label{good-defi-Wsecat}{\rm
If $p:E\rightarrow B$ and $p':E'\rightarrow B'$ are weakly
equivalent morphisms, then $\mbox{Wsecat}(p)=\mbox{Wsecat}(p').$}
\end{proposition}

\begin{proof}
We can suppose, without losing generality, that there is a
commutative square
$$\xymatrix{
{E} \ar[d]_p \ar[r]^u_{\sim } & {E'} \ar[d]^{p'} \\
{B} \ar[r]^{\sim }_{v} & {B'} }$$ Now, if $\xymatrix{B'
\ar[r]_{\sim }^{\tau '} & {F'} \ar@{>>}[r] & {0}}$ is an
$F$-factorization of the zero morphism, then an $F$-factorization
$\xymatrix{B \ar[r]_{\sim }^{\tau } & {F} \ar@{>>}[r]^w & {F'}}$
of $\tau 'v$ gives rise to $\xymatrix{B \ar[r]_{\sim }^{\tau } &
{F} \ar@{>>}[r] & {0}},$ another $F$-factorization, and a
commutative square
$$\xymatrix{
{E} \ar[d]_{\tau p} \ar[r]^u_{\sim } & {E'} \ar[d]^{\tau 'p'} \\
{F} \ar@{>>}[r]^{\sim }_{w} & {F'} }$$ Again, the result follows
considering a similar argument to that given in the proof of
Proposition \ref{salvoequivdeb}.
\end{proof}

We now see that $\mbox{Gsecat}$ and $\mbox{Wsecat}$ coincide in a
$J$-category.

\begin{theorem}\label{coinciden}{\rm
If $p:E\rightarrow B$ is any morphism, then
$$\mbox{Gsecat}(p)=\mbox{Wsecat}(p).$$}
\end{theorem}

For it we recall some useful properties about weak pullbacks.
Again we refer the reader to \cite{D}.

\begin{lemma}[\textbf{Prism Lemma for weak pullbacks}]
{\rm \cite[Prop. 2.5]{D} Consider the following diagram
$$\xymatrix{
A \ar[d]  & B \ar[d]  & C \ar[d] \\
X \ar[r] & Y \ar[r] & Z }$$ If $B\mbox{-}C\mbox{-}Z\mbox{-}Y$ is a
weak pullback, then $A\mbox{-}B\mbox{-}Y\mbox{-}X$ is a weak
pullback if and only if $A\mbox{-}C\mbox{-}Z\mbox{-}X$ is a weak
pullback.}
\end{lemma}

\begin{lemma}\label{lift-weakpullback}{\rm \cite[Lemma 3.5]{D}
 Consider a weak pullback
$$\xymatrix{
{D} \ar[d]_{g}^{\hspace{15pt}\mbox{h.p.b.}}  \ar[r] & {\bullet }
\ar@{>>}[dr] & {C} \ar[d]^{g'} \ar[l]^{\sim } \\ A \ar[rr]_f & &
{B} }$$ \noindent and let $h:X\rightarrow A$ be any morphism. Then
$h$ admits a weak lifting along $g$ if and only if $fh$ admits a
weak lifting along $g'.$}
\end{lemma}

And now the Join Theorem. This result strongly relies on the cube
axiom (J5 axiom) and therefore it does not admit a dual version.

\begin{lemma}[\textbf{Join Theorem}]\label{jointheorem}
{\rm \cite[Th. 2.7]{D} Consider the weak pullbacks
$$\xymatrix{ A \ar[d]_f^{\hspace{15pt}\mbox{h.p.b.}}  \ar[r] & X
\ar@{>>}[dr]_p & {A'} \ar[d] \ar[l]^{\sim } & & C
\ar[d]_g^{\hspace{15pt}\mbox{h.p.b.}}  \ar[r] & Y \ar@{>>}[dr]_q
& {C'} \ar[d] \ar[l]^{\sim } \\
B \ar[rr]_b & & {B'} & & B \ar[rr]_b & & {B'}}
$$
Then there is a weak pullback
$$\xymatrix{
{A*_B C} \ar[d]^{\hspace{15pt}\mbox{h.p.b.}}  \ar[r] & {\bullet }
\ar@{>>}[dr] & {A'*_{B'}C'} \ar[d] \ar[l]^(.6){\sim } \\
B \ar[rr]_b & & {B'} }$$}
\end{lemma}

\begin{proof}[Proof of Theorem \ref{coinciden}]
First suppose that $B$ is e-fibrant. We will see by induction on
$n\ge 0$ that for any map $p{\colon}E\rightarrow B$, there is weak
pullback:
$$\xymatrix{
{*^n_B E} \ar[d]_{h_n}^{\hspace{15pt}\mbox{h.p.b.}}  \ar[r] &
{\bullet } \ar@{>>}[dr] & {T^n(p)} \ar[d]^{j_n} \ar[l]^{\sim } \\
B \ar[rr]_{\Delta _{n+1}} & & {B^{n+1}} }$$ For $n=0$ it is
trivially true. Suppose the statement true for $n-1$ and consider
the diagram
$$
\xymatrix{ {*^{n-1}_B E} \ar[d]_{h_{n-1}}
\ar@{}[dr]|-(.4){\textcircled{{\scriptsize 1}}} &
{T^{n-1}(p)\times B }\ar[r]^{pr}
\ar[d]|-{j_{n-1}\times id_B} & {T^{n-1}(p)} \ar[d]_{j_{n-1}} \\
{B} \ar[r]_{\Delta _{n+1}} & {B^n\times B} \ar@{>>}[r]_{pr} &
{B^n} }
$$
where the right square is a pullback in which $pr:B^n\times
B\twoheadrightarrow B^n$ is a fibration (observe that $B$ is
e-fibrant and use (J2) axiom). Therefore this pullback is also a
homotopy pullback and a weak pullback. Now, applying the Prism
Lemma together with the induction hypothesis we deduce that
diagram {\textcircled{{\scriptsize 1}}} is also a weak pullback.

The same argument applied to the diagram
$$
\xymatrix{ {E} \ar[d]_{p} \ar[rr]^(.4){(p,p,...,p,id_E)}
\ar@{}[drr]|-{\textcircled{{\scriptsize 2}}} & & {B^n\times
E} \ar[d]|-{id_{B^n}\times p} \ar@{>>}[rr]^{pr} & & {E} \ar[d]_p  \\
{B} \ar[rr]_{\Delta _{n+1}} & & {B^n\times B} \ar@{>>}[rr]_{pr} &
& B }
$$
implies that \textcircled{{\scriptsize 2}} is a weak pullback. We
obtain the expected result by applying the Join Theorem  to the
weak pullbacks \textcircled{{\scriptsize 1}} and
\textcircled{{\scriptsize 2}}. The theorem easily follows now from
this fact together with Lemma \ref{lift-weakpullback}.

When $B$ is not e-fibrant, consider $\xymatrix{B \ar[r]_{\sim
}^{\tau } & {F} \ar@{>>}[r] & {0}}$ an $F$-factorization. Then we
have that $\mbox{Gsecat}(p)=\mbox{Gsecat}(\tau p)$ by Proposition
\ref{salvoequivdeb}. But we have already proved that
$\mbox{Gsecat}(\tau p)=\mbox{Wsecat}(\tau p)$
(=$\mbox{Wsecat}(p)$).
\end{proof}

\begin{remark}
When our category $\mathcal{C}$ does not satisfy the cube axiom
(J5), the most we can say is that $\mbox{Wsecat}(p)\leq
\mbox{Gsecat}(p)$. Indeed, a similar argument that the one used in
Theorem \ref{coinciden} using Lemma \ref{buenadefi} instead of
Lemma \ref{jointheorem}, proves that for each $n\geq 0,$ $\Delta
_{n+1}h_n$ admits a weak lifting along $j_n,$ i.e., there is a
commutative diagram
$$\xymatrix{
{*^n_B E} \ar[d]_{h_n}^{\hspace{15pt}}  \ar[r] &
{\bullet } \ar@{>>}[dr] & {T^n(p)} \ar[d]^{j_n} \ar[l]^{\sim } \\
B \ar[rr]_{\Delta _{n+1}} & & {B^{n+1}} }$$ The general case, in
which $B$ is not necessarily e-fibrant, follows easily. Now, if
$id_B$ admits a weak lifting along $h_n,$ then it is easy to check
that $\Delta _{n+1}=\Delta _{n+1}id_B$ admits a weak lifting along
$\Delta _{n+1}h_n.$ Using Lemma \ref{transitivity} below we obtain
that $\Delta _{n+1}$ admits a weak lifting along $j_n.$
\end{remark}

From now on we will denote by $\secat (p)$ both equivalent
invariants and call it the {\em sectional category} of $p.$

\section{Main properties of the sectional category}

We begin by observing that the Lusternik Schnirelmann category of
an object $B$ in $\mathcal{C}$ is the sectional category of the
zero morphism $0\rightarrow B.$ Indeed  (see \cite{D}) the $n$-th
Ganea map $p^n\colon G^nB\rightarrow B$  is precisely the $n$-th
join over $B$, $h_n:*^n_BE\rightarrow B$, of $0\rightarrow B$ and
therefore,
$$\mbox{cat}(B)=\mbox{secat}(0\rightarrow B).$$

On the other hand, given $b:B\rightarrow B'$ any morphism, we
define $\mbox{cat}(b)$ as the least integer $n\leq \infty $ such
that $b$ admits a weak lifting along $p'^n:G^nB'\rightarrow B'$.
Compare the next result with \cite{K}.

\begin{theorem}\label{superlema}{\rm
Let $p:E\rightarrow B$, $p':E'\rightarrow B'$ and $b\colon B\to
B'$ be morphisms in $\mathcal{C}$ defining a weak pullback. Then,
$$\mbox{secat}(p)\leq \min \{\mbox{cat}(b),\mbox{secat}(p')\}.$$ }
\end{theorem}
For its proof we shall need the following lemma.

\begin{lemma}\label{transitivity}{\rm \cite[Lemma 3.4]{D}
Let $f:A\rightarrow B,$ $g:C\rightarrow B$ and $h:D\rightarrow B$
be morphisms. If $f$ admits a weak lifting along $g$ and $g$
admits a weak lifting along $h$, then $f$ admits a weak lifting
along $h.$}
\end{lemma}

\begin{proof}[Proof of \ref{superlema}]
By induction, using repeatedly the Join Theorem (Lemma
\ref{jointheorem}) on the given weak pullback

$$\xymatrix{
{E} \ar[d]_{p} \ar[r] & {\bullet }
\ar@{>>}[dr] & {E'} \ar[d]^{p'} \ar[l]^(.6){\sim } \\
B \ar[rr]_b & & {B'} }$$ we obtain, for every $n\geq 0,$
 a weak pullback of the form
$$\xymatrix{
{*^n_BE} \ar[d]_{h_n} \ar[r] & {\bullet }
\ar@{>>}[dr] & {*^n_{B'}E'} \ar[d]^{h'_n} \ar[l]^(.6){\sim } \\
B \ar[rr]_b & & {B'} }$$ Hence, if $\mbox{secat}(p')\leq n$,
$h'_n$ admits a weak section:
$$\xymatrix@R=.4cm@C=.4cm{
 & & {*^n_{B'}E'} \ar[dd]^{h'_n} \ar[dl]_{\eqr{45}} \\
 & {\bullet } \ar@{>>}[dr] & \\
 {B'} \ar@{.>}[ur]^s \ar[rr]_{id_{B'}} &  & {B'} } $$
In particular, $b:B\rightarrow B'$ admits a weak lifting along
$h'_n$ through the morphism $sb:B\rightarrow {\bullet }.$ By Lemma
\ref{lift-weakpullback}, $h_n$ admits a weak section and
$\mbox{secat}(p)\leq n.$

Now suppose that $\mbox{cat}(b)\leq n,$ that is, $b$ admits a weak
lifting along $p'^n:G^nB'\rightarrow B'$. Consider the following
diagram obtained by simply choosing any $F$-factorization of $p'$:
$$\xymatrix{
{0} \ar[d] \ar[r] & {\bullet }
\ar@{>>}[dr] & {E'} \ar[d]^{p'} \ar[l]^(.6){\sim } \\
{B'} \ar[rr]_{id} & & {B'} }$$ As this is not in general a weak
pullback, apply this time Lemma \ref{buenadefi} inductively to
obtain that $p'^n:G^nB'\rightarrow B'$ admits a weak lifting along
$h'_n:*^n_{B'}E'\rightarrow B'.$ Finally, by Lema
\ref{transitivity} we conclude that $b$ admits a weak lifting
along $h'_n:*^n_{B'}E'\rightarrow B',$ which by Lemma
\ref{lift-weakpullback}, is equivalent to the fact that
$h_n:*^n_{B}E\rightarrow B$ admits a weak section.
\end{proof}

Even if our data is not a weak pullback, we can prove a similar
result. Compare with \cite{K}.

\begin{theorem}\label{lema2}{\rm
Let $p:E\rightarrow B$ and $p':E'\rightarrow B$ be morphisms in
$\mathcal{C}$. If $p$ admits a weak lifting along $p',$ then
$\mbox{secat}(p')\leq \mbox{secat}(p).$ In particular,
$$\mbox{secat}\hspace{1pt}(p)\leq
\mbox{cat}\hspace{1pt}(B).$$ Moreover, if $p:E\rightarrow B$
admits a weak lifting along the zero morphism $0\rightarrow B$ (in
particular, when $E$ is weakly contractible, i.e., $E$ and $0$
 are weakly equivalent) then
$\mbox{secat}\hspace{1pt}(p)=\mbox{cat}\hspace{1pt}(B)$.}
\end{theorem}

\begin{proof}
For the first assertion, apply  Lemma \ref{buenadefi} inductively
to the diagram
$$\xymatrix{
{E} \ar[d]_{p} \ar[r] & {\bullet }
\ar@{>>}[dr] & {E'} \ar[d]^{p'} \ar[l]^(.6){\sim } \\
B \ar[rr]_{id_B} & & {B} }$$ to conclude that, for every $n\geq
0,$ $h_n$ admits a weak lifting along $h'_n$. If $\secat(p)\le n$,
$id_B$ admits a weak lifting along $h_n$ and, by Lemma
\ref{transitivity}, $id_B$ admits a weak lifting along $h'_n$.
Hence, $\secat(p')\le n$.

On the other hand, recall that
$\mbox{cat}(B)=\mbox{secat}(0\rightarrow B)$ and observe that the
zero morphism  admits a weak lifting along any morphism. Thus,
$\mbox{secat}\hspace{1pt}(p)\leq \mbox{cat}\hspace{1pt}(B).$
Finally note that, if $E$ is a weakly trivial object, by Lemma
\ref{lifting-lemma}, $p$ admits a weak lifting along
$0:0\rightarrow B$.
\end{proof}

\subsection{Modelization functors.}

We now study the behaviour of $\secat$ through a {\em modelization
functor}. Recall from \cite{D} that a covariant functor $\mu
:\mathcal{C}\rightarrow \mathcal{D}$ between categories satisfying
(J1)-(J4) axioms is called a modelization functor if it preserves
weak equivalences, homotopy pullbacks and homotopy pushouts. We
say that $\mu $ is \textit{pointed} if  $\mu (0)=0.$ If $\mu
:\mathcal{C}\rightarrow \mathcal{D}$ is contravariant, it is said
to be a modelization functor if the corresponding covariant
functor $\mu :\mathcal{C}^{op}\rightarrow \mathcal{D}$ is a
modelization functor. Here we prove:

\begin{theorem}\label{comparison}{\rm
If $\mu :\mathcal{C}\rightarrow \mathcal{D}$ is a modelization
functor between $J$-categories, then for any morphism
$p:E\rightarrow B$ of $\mathcal{C}$
$$\mbox{secat}(\mu (p))\leq \mbox{secat}(p)$$}
\end{theorem}

For it we shall need the following

\begin{lemma}\label{tecnico}{\rm \cite[Prop. 6.7]{D}
Let $\mu :\mathcal{C}\rightarrow \mathcal{D}$ be a modelization
functor and let $j:A*_BC\rightarrow B$ denote the join map of
$f:A\rightarrow B$ and $g:C\rightarrow B$ . Then, there is a
commutative diagram
$$\xymatrix{
{\mu (A*_BC)} \ar[drr]^{\mu (j)} & & \\
{\bullet } \ar[u]_{\sim } \ar[d]^{\sim } \ar[rr] & & {\mu (B)} \\
{\mu (A)*_{\mu (B)}\mu (C)} \ar[urr]_{j'} & & }$$ \noindent where
$j'$ denotes the join morphism of $\mu (f)$ and $\mu (g).$}
\end{lemma}

\begin{proof}[Proof of Theorem \ref{comparison}]
In view of Lemma \ref{transitivity} it is sufficient to prove
that, for each $n$,  $h_n^{\mu (p)}$ admits a weak lifting along
$\mu (h_n^{p})$
$$\xymatrix@R=.4cm@C=.4cm{
 & & {\mu (*^n_BE)} \ar[dl]^{\sim } \ar[dd]^{\mu (h_n^p)} \\
& {\bullet } \ar@{>>}[dr] & \\
{*^n_{\mu (B)}\mu (E)} \ar[ur] \ar[rr]_{h_n^{\mu (p)}} & & {\mu
(B)}}$$ \noindent where $h_n^p$ and $h_n^{\mu (p)}$ are the $n$-th
join morphisms $p$ and $\mu (p)$ respectively. For $n=0$ is
trivially true. By assuming the assertion true for $n-1,$ and
choosing any $F$-factorization of $\mu (p)$ we obtain a
commutative diagram of the form
$$\xymatrix{
{*^{n-1}_{\mu (B)}\mu (E)} \ar[d] \ar[rr]^{h_{n-1}^{\mu (p)}} & &
{\mu (B)} \ar[dd]^{id} & & {\mu (E)} \ar[ll]_{\mu (p)}
\ar[d]_{\sim } \\
{\bullet } \ar@{>>}[drr] & & & & {\bullet } \ar@{>>}[dll] \\
{\mu (*^{n-1}_{B}E)} \ar[u]_{\sim } \ar[rr]_{\mu (h_{n-1}^p)} & &
{\mu (B)} & & {\mu (E)} \ar[ll]^{\mu (p)} \ar[u]^{\sim } }$$ By
Lemma \ref{buenadefi} $h_n^{\mu (p)}$ admits a weak section along
the join morphism of $\mu (h_{n-1}^p)$ and $\mu (p):$
$$\xymatrix{
{*^n_{\mu (B)}\mu (E)} \ar[d]_{h_n^{\mu (p)}} \ar[r] & {\bullet }
\ar@{>>}[dr] & {\mu (*^{n-1}_{B}E)*_{\mu (B)}\mu (E)} \ar[d] \ar[l]^(.75){\sim } \\
{\mu (B)} \ar[rr]_{id} & & {\mu (B)} }\eqno(3)$$ On the other
hand, applying Proposition \ref{tecnico} above to the morphisms
$h_{n-1}^p:*^{n-1}_BE\rightarrow B$ and $p:E\rightarrow B,$ we
obtain a commutative diagram
$$\xymatrix{
{\bullet } \ar[rr]_{\sim } \ar[d]^{\sim } & & {\mu (*^n_BE)}
\ar[d]^{\mu (h_n^p)} \\
{\mu (*^{n-1}_{B}E)*_{\mu (B)}\mu (E)}  \ar[rr] & & {\mu (B)}}$$
Taking any $F$-factorization of $\mu (h_n^p)$ and applying Lemma
\ref{lifting-lemma} we deduce that the join morphism $\mu
(*^{n-1}_{B}E)*_{\mu (B)}\mu (E)\rightarrow \mu (B)$ admits a weak
lifting along $\mu (h_n^p)$. Finally, by Lemma \ref{transitivity}
applied to $(3)$, we conclude the inductive step.
\end{proof}

\begin{remark}
Observe that, for the proof of Theorem \ref{comparison} we have
used the Ganea-type version of sectional category. If (J5) axiom
is not satisfied, then using similar arguments we can also obtain
the same result for the Whitehead-type version of sectional
category. The same also applies for the remaining results of this
section.
\end{remark}

\begin{corollary}\label{igualdad}{\rm
Consider $\mu :\mathcal{C}\rightarrow \mathcal{D}$ and $\nu
:\mathcal{D}\rightarrow \mathcal{C}$ modelization functors between
$J$-categories and let $p:E\rightarrow B$ be a morphism in
$\mathcal{C}$ such that $\nu (\mu (p))$ is weakly equivalent to
$p$. Then
$$\mbox{secat}(\mu (p))=\mbox{secat}(p)$$}
\end{corollary}

As an example we apply the theorem above to the {\em abstract
topological complexity} of a given object. For any e-fibrant
object $B$ we define its topological complexity, $\mbox{TC}(B)$ as
the sectional category of the diagonal morphism $\Delta
_B:B\rightarrow B\times B$. If $B$ is not e-fibrant consider any
$F$-factorization $\xymatrix{B \ar[r]^{\sim } & {F} \ar@{>>}[r] &
{0}}$ and set
$$\mbox{TC}(B):=\mbox{TC}(F).$$ Then $\mbox{TC}(B)$ does not depend on the
e-fibrant object $F;$ indeed, if we take another
$F$-fac\-tor\-ization $\xymatrix{B \ar[r]^{\sim } & {F'}
\ar@{>>}[r] & {0}},$ then there exists a weak equivalence
$h:F'\stackrel{\sim }{\rightarrow }F$ (see the proof of Lemma
\ref{good-defi-whit}). The naturality of the diagonal morphism
applied to $h$ together with the fact that $h\times h:F'\times
F'\stackrel{\sim }{\rightarrow }F\times F$ is a weak equivalence
(by the dual of the Gluing Lemma \cite[II.1.2]{B}) prove that
$\Delta _{F}:F\rightarrow F\times F$ and $\Delta
_{F'}:F'\rightarrow F'\times F'$ are weakly equivalent morphisms.
Therefore $$\mbox{TC}(F)=\mbox{secat}(\Delta
_{F})=\mbox{secat}(\Delta _{F'})=\mbox{TC}(F').$$ $\mbox{TC}(B)$
neither depends on the weak type of $B;$ given $f:B\stackrel{\sim
}{\rightarrow }B'$ a weak equivalence, if we consider an
$F$-factorization $\xymatrix{{B'} \ar[r]_{\sim }^{\tau '} & {F'}
\ar@{>>}[r] & {0}}$, then any $F$-factorization of the composite
$\tau 'f:B\rightarrow F'$
$$\xymatrix{ {B} \ar[rr]^{\tau 'f}_{\sim } \ar[dr]^{\sim }_{\tau } & & {F'}
\\ & {F} \ar@{>>}[ur]_g & }$$ \noindent gives rise to a trivial
fibration $g:F\stackrel{\sim }{\twoheadrightarrow }F',$ which
shows that
$$\mbox{TC}(B)=\mbox{TC}(F)=\mbox{TC}(F')=\mbox{TC}(B').$$

\begin{theorem}{\rm
For any pointed modelization functor  $\mu :\mathcal{C}\rightarrow
\mathcal{D}$ and any object $B$, $$\mbox{TC}(\mu (B))\leq
\mbox{TC}(B)$$}
\end{theorem}

\begin{proof}
Taking into account that $\mu $ preserves weak equivalences and
$\mbox{TC}$ does not depend on the weak type, we can suppose
without losing generality that $B$ is an e-fibrant object. Since
$\mu (B)$ need not be e-fibrant we consider any $F$-factorization
$$\xymatrix{{\mu (B)} \ar[rr] \ar[dr]^{\sim }_{\tau } & & {0}
\\ & {F} \ar@{>>}[ur] & }$$
\noindent so that $\mbox{TC}(\mu (B))=\mbox{TC}(F).$ Now take the
following commutative cube:
$$\xymatrix@R=.4cm@C=.4cm{
  & {\mu (B\times B)} \ar[rrr]^{\mu (pr_2)} \ar@{.>}'[dd]^{\,\omega }[dddd] \ar[ddl]_{\mu (pr_1)}
      &  &  & {\mu (B)} \ar[dddd]^{\tau }_{\sim} \ar[ddl]     \\
                &  &  & \\
  {\mu (B)} \ar[rrr] \ar[dddd]_{\tau }^{\sim }
      &  &  & {0} \ar[dddd]^(.3){id}_(.3){\sim } \\
                &  &  &   \\
  & {F\times F} \ar@{>>}'[rr]_(.7){pr_2}[rrr]  \ar@{>>}[ddl]_{pr_1} &  &  &
  {F} \ar@{>>}[ddl]
       \\
                &  &  &\\
  {F} \ar@{>>}[rrr] &  &  &  {0} \\
}$$ \noindent where $pr_1$ and $pr_2$ denote the projection
morphisms.  As $\mu$ is a pointed modelization functor, the top
face is a homotopy pullback. On the other hand, the bottom face is
a strict pullback (and a homotopy pullback) and $\omega =(\tau \mu
(pr_1),\tau \mu (pr_2))$ is the induced morphism from the
universal property of the pullback. Since the top and bottom faces
are homotopy pullbacks and the unbroken vertical morphisms are
weak equivalences, by \cite[Cor. 1.12]{D} (or the dual of the
Gluing Lemma \cite[II.1.2]{B}) we have that $\omega $ is also a
weak equivalence. From the following commutative diagram
$$\xymatrix{
{\mu (B)} \ar[rr]^{\tau }_{\sim } \ar[d]_{\mu (\Delta _B)} & &
{F} \ar[d]^{\Delta _{F}} \\
{\mu (B\times B)} \ar[rr]^{\sim }_{\omega } & & {F\times F} }$$
\noindent we deduce that $\mu (\Delta _B)$ and $\Delta _{F}$ are
weakly equivalent morphisms. Then, by Proposition
\ref{salvoequivdeb} we have that $\mbox{TC}(\mu
(B))=\mbox{secat}(\Delta _{F}) =\mbox{secat}(\mu (\Delta _B))$
while, by Theorem \ref{comparison}, $\mbox{secat}(\mu (\Delta
_B))\leq \mbox{secat}(\Delta _B)=\mbox{TC}(B)$.
\end{proof}

\begin{corollary}\label{igualdad2}{\rm
Consider $\mu :\mathcal{C}\rightarrow \mathcal{D}$ and $\nu
:\mathcal{D}\rightarrow \mathcal{C}$ pointed modelization functors
and let $B$ be an object in $\mathcal{C}$ such that $\nu (\mu
(B))$ is weakly equivalent to $B.$ Then
$$\mbox{TC}(\mu (B))=\mbox{TC}(B)$$}
\end{corollary}

\section{Some applications.}

We start by an immediate application in rational homotoy theory. A
classical fact \cite[\S8]{bkan} assures the existence of an
adjunction
$$\xymatrix{ {\text{\bf CDGA}^\varepsilon}&{\text{\bf SSet}^*}
\ar@<1ex>[l]^{\langle\,\cdot\,\rangle} \ar@<1ex>[l];[]^{A_{PL}}&
\\}$$ between the
categories of augmented commutative differential graded algebras
over a field $\mathbb{K}$ of characteristic zero, and pointed
simplicial sets. The category ${\text{\bf SSet}^*}$ is known to be
a $J$-category endowed with Kan fibrations, injective maps and
maps realizing to homotopy equivalences \cite[Chap.III\S.3]{Q},
\cite[Prop.A.8]{D}. The category ${\text{\bf CDGA}^\varepsilon}$
is also a (proper) closed model category \cite[\S4]{bkan} (and
thus J1-J4 are satisfied) in which fibrations are surjective
morphisms, weak equivalences are morphisms inducing homology
isomorphisms (the so called ``quasi-isomorphisms") and
cofibrations are ``relative Sullivan algebras" \cite[\S14]{fht},
i.e., inclusions $A\to A\otimes\Lambda V$ in which $\Lambda V$
denotes the free commutative algebra generated by the graded
vector space $V$ and the differential on $A\otimes \Lambda V$
satisfies a certain ``minimality" condition. However, this is NOT
a $J$-category and the Eckmann-Hilton dual of a partial version of
the cube axiom is satisfied when restricting to $1$-connected
algebras \cite[A.18]{D}. The functors ${\langle\,\cdot\,\rangle}$
and $A_{PL}$ do not in general respect weak equivalences although
${\langle\,\cdot\,\rangle}$ sends cofibrations to fibrations and
${\langle\,\cdot\,\rangle}$ can be slightly modified to send
fibrations to cofibrations \cite[\S8]{bkan}.  Therefore, as they
stand, they are not modelization functors. However, it is also
known \cite[\S8,9]{bkan} that, restricting those functors to the
categories
$$\xymatrix{ {\text{\bf CDGA}^1_{cf\Q}}&{\text{\bf Kan-Complexes}^1_\Q}\ar@<1ex>[l]
\ar@<1ex>[l];&\\}$$ of cofibrant $1$-connected commutative
differential graded algebras of finite type over $\Q$ (known as
Sullivan algebras \cite[\S12]{fht}) and $1$-connected rational Kan
complexes of finite type, then they do preserve weak equivalences
and via \cite[Prop.6.5]{D} they are modelization functors.

On the other hand, in \cite[Ch.8]{fa}, Fass\`o introduce, for a
map of finite type $1$-connected CW-complexes, or equivalently for
a simplicial map of finite type $1$-connected Kan complexes
$E\stackrel{p}{\to}B$ the {\em rational sectional category} of
$p$, $\secat_0(p)$ which can be seen as the sectional category in
the opposite category of $\text{\bf CDGA}^1_{cf\mathbb{Q}}$ of
$A_{PL}( p_{\mathbb{Q}}),$ being $p_{\mathbb{Q}}$ the map in
${\text{\bf Kan-Complexes}^1_\Q}$ obtained by rationalization
\cite[\S11]{bkan}. Thus,  by Corollary \ref{igualdad},
$$\secat_0(p)=\secat(p_{\mathbb{Q}})$$

\bigskip
Our second application concerns localization functors. Let $P$ be
a (possibly empty) set of primes and
$$
(-)_P\colon \text{\bf CW}_{\mathcal{N}}\longrightarrow \text{\bf
CW}_{\mathcal{N}}
$$
denotes the $P$-localization functor (see \cite[\S2]{hmr} or
\cite[Chap.III]{ar} where it is shown that localization can chosen
to be a functor as it stands, not just in the homotopy category)
in the pointed category of spaces of the homotopy type of
nilpotent CW-complexes. Then, this functor sends homotopy pushouts
to homotopy pushouts and homotopy pullbacks (if the chosen
homotopy pullback stays in this category) to homotopy pullbacks
\cite[\S7]{hmr}. (Note that, considering closed cofibrations,
Hurewicz fibrations and homotopy equivalences, the category of
well pointed topological spaces ${\text{\bf Top}^*}$ has the
structure of a $J$-category; see \cite[Thm.11]{S3} for axioms
(J1)-(J4) plus \cite[Thm.25]{M} for (J5)). Thus, even though
strictly speaking this is not a modelization functor as it is
defined on a certain subcategory of ${\text{\bf Top}^*},$ the
arguments in Theorem \ref{comparison} could be followed mutatis
mutandi as long as all constructions there remain within our
category. But this is in fact the case as the homotopy pullback
(or pushout) of two maps in $\text{\bf CW}_{\mathcal{N}}$ can be
chosen to live also in this category \cite[\S7]{hmr}. Hence,
$$
\secat(f_P)\le \secat f.
$$
However, the situation is drastically different in the general
case as all sort of possible $P$-localizations (extending the one
on nilpotent complexes) do not, in general, preserve homotopy
pullbacks and homotopy pushouts.

Here, we consider the Casacuberta-Peschke localization functor on
${\text{\bf Top}^*}$ \cite{cp} and start by setting some notation.
Given a group $G$ we denote by P[G]  the ring localization of the
group ring $\mathbb{Z}_PG$ obtained by inverting all of the
elements $1 + g +\cdots+ g^{n-1}$, where $g \in G$, and $(n,p)=1$
for any $p\in P$ (see \cite[\S2]{cp}).

Following \cite{p} we say that a $P$-torsion group $G$ is an {\em
acting group} for a space $X$ if there is an epimorphism $f\colon
\pi_1X\twoheadrightarrow G$ such that, for each $m\ge 2$, the
action $\pi_1X\to Aut(\pi_mX)$ factors through $G$.

\begin{proposition}{\rm
Let $f\colon X\to Y$ a map for which:
\begin{enumerate}
\item[(i)] $\pi_1(*^n_Yf)\colon
\pi_1(*_Y^nX)\stackrel{\cong}{\longrightarrow} \pi_1Y$ is an
  isomorphism of $P$-local groups for any $n\ge 0$.

\item[(ii)] $\pi_1(*^n_YX)$ and $\pi_1Y$ have a common acting group
$G$ for any $n\ge 0$.

\item[(iii)] If we denote $\pi_1Y$ by $\pi$, the morphism
$\mathbb{Z}_P\pi\to P[\pi]$ induce isomorphisms on homology with
local coefficients $H_*(-;\mathbb{Z}_P\pi)\to H_*(-;P[\pi])$.
\end{enumerate} Then,
$$\secat(f_P)\le \secat (f).$$}
\end{proposition}

\begin{proof}
Again, note that the argument in Theorem \ref{comparison} could be
applied if, for any $n\ge 1$, there is a homotopy commutative
diagram of the form:
$$
\xymatrix{{(*^{n-1}_YX)_P*_{Y_P}X_P}\ar[rd]_{{(h_{n-1})_P*f_P}}\ar[rr]^{\simeq}
&&{(*^n_YX)_P}\ar[ld]^{(h_n)_P}\\
&Y_P&}
$$
To this end, an inductive process,  as in \cite[Prop.6.7]{D} will
work as long as the following two conditions hold:

(1) The localization of the homotopy pullback
$$\xymatrix{Q_n\ar[r]\ar[d]&X\ar[d]^{f}\\{*^{n-1}_YX}\ar[r]_(.7){h_{n-1}}&Y}$$
$$\xymatrix{(Q_n)_P\ar[r]\ar[d]&X_P\ar[d]^{f_P}\\{(*^{n-1}_YX)_P}\ar[r]_(.7){{h_{n-1}}_P}&Y_P}$$

is again a homotopy pullback.

(2) The localization of the homotopy pushout
$$\xymatrix{Q_n\ar[r]\ar[d]&X\ar[d]\\{*^{n-1}_YX}\ar[r]&{*^n_YX}}$$
$$\xymatrix{(Q_n)_P\ar[r]\ar[d]&X_P\ar[d]\\{(*^{n-1}_YX)_P}\ar[r]&{(*^n_YX)_P}}$$
is again a homotopy pushout.

However, by hypothesis, we may apply \cite[Thm.4.3]{p} to prove
statement (1) (res\-pec. \cite[Thm.2.1]{p} to prove (2)).
\end{proof}

\end{document}